\documentclass[10pt,twoside]{article}
\usepackage{amssymb}
\usepackage[mathscr]{eucal}
\usepackage{eufrak}
\usepackage{amsmath}
\usepackage{mathrsfs}

\def\O{\Omega}

\def\R{\Bbb R}
\def\N{\Bbb N}
\def\o{\"{o}}
\def\à{\`{a}}
\def\è{\`{e}}
\def\ì{\`{i}}
\def\ù{\`{u}}
\def\ò{\`{o}}
\def\é{\'{e}}

\def\dy{\displaystyle}
\def\ve{\varepsilon}

\def\be{\begin{equation}}
\def\ba{\begin{array}}
\def\ea{\end{array}}
\def\ee{\end{equation}}
\def\vs1{\vspace{1ex}}
\def\vp{\varphi}
\def\ov{\overline}
\def\po{\partial\Omega}
\font\sc=cmcsc10
\title{Global $L^r$-estimates and regularizing effect 
 for solutions to the  $p(t,x)\,$-Laplacian  systems}\author{\sc F. Crispo, P. Maremonti\thanks{
Dipartimento di Matematica e Fisica, Universit\`{a} degli Studi della Campania ``Luigi Vanvitelli'', viale Lincoln 5, 81100 Caserta,
 Italy. 
 },  M. R\r u\v{z}i\v{c}ka\thanks{Mathematisches Institut, Universit\"{a}t Freiburg, Ernst-Zermelo-Str.~1, D-79104 Freiburg i. Br., Germany.}} 
 \date{}
\begin{document}
\maketitle
\noindent{\bf Abstract}  - {\small We consider the initial boundary value
  problem for the $p(t,x)$-Laplacian system in a bounded domain
  $\Omega$. If the initial data belongs to $L^{r_0}$, $r_0\geq 2$, we give a global $L^{r_0}(\O)$-regularity result uniformly in
  $t>0$ that, in the particular case ${r_0}=\infty$, implies a maximum
  modulus theorem.  Under the assumption $p_-=\inf p(t,x)>\frac{2n}{n+r_0}$, we also state 
   $L^{r_0}-L^{r}$ estimates for the solution, for $r\geq r_0$.
Complete proofs of the results presented here are given in the paper \cite{CMR}.}
 \vskip -0.7true cm\noindent
\newcommand{\red}{\protect\bf}
\renewcommand\refname{\centerline
{\red {\normalsize \bf References}}}
\newtheorem{ass}
{\bf Assumption}[section]
\newtheorem{defi}
{\bf Definition}[section]
\newtheorem{tho}
{\bf Theorem}[section]
\newtheorem{rem}
{\sc Remark}[section]
\newtheorem{lemma}
{\bf Lemma}[section]
\newtheorem{coro}
{\bf Corollary}[section]
\newtheorem{prop}
{\bf Proposition}[section]
\renewcommand{\theequation}{\thesection .
\arabic{equation}}
\setcounter{section}{1}
\section*{\normalsize 1. Introduction}
\renewcommand{\theequation}{1.\arabic{equation}}
\renewcommand{\thetho}{1.\arabic{tho}}
\renewcommand{\thedefi}{1.\arabic{defi}}
\renewcommand{\therem}{1.\arabic{rem}}
\renewcommand{\theprop}{1.\arabic{prop}}
\renewcommand{\thelemma}{1.\arabic{lemma}} \setcounter{equation}{0}
\setcounter{lemma}{0} \setcounter{defi}{0} \setcounter{prop}{0}
\setcounter{rem}{0} \setcounter{tho}{0} 

The aim of this note is to present a global $L^{r_0}$-regularity
result, uniformly in $t$, and, to study the regularizing effect for
solutions of the following $p(t,x)$-Laplacian system \be\label{PF}
\begin{array}{rll}
\vs1\dy u_t-\nabla\cdot  (|\nabla u|^{p(t,x)-2}\nabla u)&\!\!\!\!= 0,&\hskip-0.1cm\textrm{ in
}(0,\infty)\times\O,
\\
  \dy\vs1u(t,x)&\!\!\!\!=0,&\hskip-0.1cm \textrm{ on }
(0,\infty)\times\partial\O,\\\dy u(0,x)&\!\!\!\!=u_\circ(x),&\hskip-0.1cm\mbox{ in
}\O.
\end{array}\ee 
Here  $\O\subset\R^n$, $n\geq 2$, is a bounded domain, 
$u:(0,\infty)\times\O\to \R^N,\,N\geq1$ is a scalar or vector field, $u_t$
is the derivative of $u$ with respect to time, and the exponent
$p:=p(t,x)$ is bounded
\be\label{ass1}\ba{ll}\dy\vs1 1<p_{-}:=\inf p(t,x)\leq p_-(t):=\inf_{\O} p(t,x)\leq p(t,x)\\
\dy \hskip4cm\leq p_+(t):=\sup_{\O} p(t,x)\leq p_+:=\sup p(t,x)< \infty\,,\ea\ee
 and log-H\o lder continuous in $(t,x)$, which means that 
 there exists a constant $c_1$ such that
 \be\label{ass2}|p(t,x)-p(\tau,y)|\leq \frac{c_1}{\log(e+\frac {1}{|t-\tau|+|x-y|})}\ee
is satisfied for all $(t,x)$ and $(\tau,y)$  in $(0,\infty)\times\O$.\par
\par 
 System \eqref{PF} belongs to the class of
partial differential equations with non-standard growth.  One of the
principal reason of its interest is the connection with the more
intricate system modeling the motion of electro-rheological fluids
with shear dependent viscosity, introduced in \cite{RR,RR1}. For
results on these fluids in the evolutionary case we refer to
\cite{CrispoP,DR, nagele,Ruz1,Ruz2}.\par
The well-posedness of problem \eqref{PF} is rather recent. The first
existence result, which allows to study this system by using classical
tools of monotone operators theory, appeared in
\cite{AS}. Subsequently in \cite{DNR} and in the same period in
\cite{AZ1}, the existence was extended by removing the lower bound
$p_-> \frac{2n}{n+2}$ and, in \cite{AZ1}, also the log-H\o lder
continuity condition.  The technique of these two papers are
completely different and the solutions obtained seems to be not
comparable, unless a log-H\o lder continuity condition is assumed, as
proved in \cite{AZ2}.

\par The following main theorem shows a global $L^{r_0}$-regularity
result, uniformly in $t$, for a weak solution (for definition see next section).
\begin{tho}\label{mmtp}{\sl Let $p$ satisfy assumptions \eqref{ass1}--\eqref{ass2} and let $r_0\in [2,\infty]$. Then, for all
    $u_\circ\in L^{r_0}(\O)$ the following estimate holds for the unique weak
    solution $u$ of \eqref{PF} \be\label{mmtp1}\|u(t)\|_{r_0}\leq
    \|u_\circ\|_{r_0}, \qquad \forall t\geq 0.\ee
    Moreover, for all $T\in(0,\infty)$,
    $u\in C([0,T]; L^{r}(\O))$, for $r=2$ and for any $r\in [2,r_0)$ if
    $r_0\in (2,\infty]$. Finally, if
    $r_0\in [2,\infty)$,  then $\dy\lim_{t\to 0}\|u(t)-u_\circ\|_{r_0}=0$.}
\end{tho}
In the case ${r_0}=\infty$, analogous results are proved, with a constant
exponent, in \cite{DB,DBH, DBUV} for equations, in \cite{choe} locally
for systems and, finally, globally for systems in \cite{CMADE}. Again
for equations but with a variable exponent $p(t,x)$, a similar result
is proved in \cite{AS}, where a more general non-linear parabolic
equation is considered. Finally, still for the $p(t,x)$-system, the
result is proved locally in the recent paper \cite{BD}. \par In \cite{CMR} we prove
Theorem \ref{mmtp} in a direct way, via a duality technique that makes
use of a suitable, still quasi-linear, ``adjoint problem''. This means
that we reach our result without the investigation on high regularity
properties of solutions, that could ensure boundedness by
embedding. On the other hand this last question is still an
interesting open problem. Actually in the framework of parabolic equations with non constant exponents, the unique regularity results known
are the one in the pioneering paper \cite{AMS}, where a partial H\o
lder continuity of the spatial gradient is proved, and the result in
\cite{BD}, where a local H\o lder continuity of the spatial gradient
is proved. The results in \cite{AMS} and \cite{BD} are proved under the assumption, quite
natural in the regularity theory, of $p_-> \frac{2n}{n+2}$.  \par Note
that we employ the duality as we work on weak solutions, so we do not have enough regularity
to give sense to an expression of the kind $<u_t, u|u|^{r-2}>$,
$r\in [2,\infty)$. Once Theorem \ref{mmtp} is obtained, the previous
quantity makes always sense if the initial data is in $L^\infty(\O)$,
and then for any initial data throughout a limit procedure.  
\par Our next result is Theorem \ref{mmtpx} where we get
informations on the behavior of the $L^{r}(\O)$-norm of the solution, corresponding to an initial data in $L^{r_0}(\O)$, $2\leq r_0\leq r\leq \infty$, in a right neighborhood of $t=0$ as well as for $t\to \infty$. This kind of property is known as regularizing effect.  
\begin{tho}\label{mmtpx}{\sl Let $u_\circ\in L^{r_0}(\O)$ with $r_0\in [2,\infty)$,
    and let $p$ satisfy assumptions \eqref{ass1}--\eqref{ass2} with 
$p_->\frac{2n}{n+r_0}$. Then, for all  $r\in [r_0,\infty)$, there exist nonnegative constants $\gamma_-(r)$, $\gamma_+(r)$ such that the weak solution
    $u$ of \eqref{PF} satisfies the following estimate
  \be\label{mmtp2}
\|u(t)\|_{r}\!\leq c(\|u_\circ\|_{r_0})(t^{-\gamma_-(r)}+t^{-\gamma_+(r)}),\ \forall t>0.
     \ee 
     If $p_->n$ then $r=\infty$ is allowed. }
\end{tho}
\begin{rem} {\rm In the proof of Theorem \ref{mmtpx} in \cite{CMR}  we obtain a precise expression of $\gamma_-(r)$ and $\gamma_+(r)$ if $p_-\geq n$. Actually for $r\in [r_0,\infty)$ 
$\gamma_-(r):=\frac{(r-r_0)n}{r(r_0p_--2n+np_-)}\, $ and $\,\gamma_+(r):=\frac{(r-r_0)np_-}{rp_+(r_0p_-\!\!-2n+np_-)}$. 
If $p_->n$ then $r=\infty$ is allowed and 
  $\gamma_-:=\frac{n}{r_0p_--2n+np_-}\, $ and $\,\gamma_+:=\frac{np_-}{p_+(r_0p_--2n+np_-)}$.  Note that these exponents completely agree with those given in \cite{HV,Ver,MP1,MP2} for the constant exponent case.  
  \par 
    If $p_-<n$ we could also give the expression of $\gamma_-(r)$ and  $\gamma_+(r)$. On the other hand, as we treat this case by an iterative argument, it would become very tedious. }
\end{rem}
\begin{coro}\label{nc1} 
{\sl 
Let the assumptions of Theorem \ref{mmtpx} be satisfied. Then $u\in C([0,T]; L^{r_0}(\O))$. Further 
\be\label{bd}
\|u(t)\|_{r_0}\leq \|u(s)\|_{r_0}, \ \forall t >s\geq 0.\ee}    
\end{coro}
\begin{rem}{ \rm 
Having the result of Corollary \ref{nc1} at disposal, we can state estimate \eqref{mmtp2} in the form of parabolic semigroup as
$$\|u(t)\|_{r}\leq c(\|u(s)\|_{r_0}) ( (t-s)^{\gamma_-(r)}+  (t-s)^{\gamma_+(r)}),\ \forall t>s>0.$$}
\end{rem}

\section*{\normalsize 2.
Ideas of the proofs}
\renewcommand{\theequation}{2.\arabic{equation}}
\renewcommand{\thetho}{2.\arabic{tho}}
\renewcommand{\thedefi}{2.\arabic{defi}}
\renewcommand{\therem}{2.\arabic{rem}}
\renewcommand{\theprop}{2.\arabic{prop}}
\renewcommand{\thelemma}{2.\arabic{lemma}}
\setcounter{equation}{0} \setcounter{lemma}{0} \setcounter{defi}{0}
\setcounter{prop}{0} 
\setcounter{rem}{0} \setcounter{tho}{0} 
Let us introduce few notation related to function spaces. We use standard notation for Sobolev and Bochner spaces. For spaces with variable exponents 
we mainly refer to the notation in \cite{AS} and \cite{DHHR}. 
Given a bounded and log-H\o lder continuous exponent $q\in \Omega_T$, with $q_->1$, for any $t>0$ we introduce the Banach space
$$V^q_t(\O):=\{u\in L^2(\O)\cap W_0^{1,1}(\O), \nabla u\in L^{q(t,\cdot)}(\O)\}\,,$$
with norm 
$$\|u\|_{V^q_t(\O)}:=\|u\|_{2,\O}+\|\nabla u\|_{q(t,\cdot), \O},$$
and we denote by $(V^q_t)'$ its dual. 
We further introduce, for all $T\in (0,\infty)$, the Banach space  
$$
  X^q(\O_T):=\{u\in L^2(\O_T), u (t,\cdot)\in V^q_t(\O)
  \mbox{ a.e. }t\in [0,T], \nabla u\in L^{q(\cdot)}(\O_T)\},
$$
with norm
$$\|u\|_{X^q(\O_T)}:=\|u\|_{2, \O_T}+\|\nabla u\|_{q(\cdot), \O_T}\,,$$
and we denote by $(X^q(\O_T))'$ its dual. Last we introduce the space
 $$W^q(\O_T):=\{u\in X^q(\O_T), u_t\in (X^q(\O_T))'\}.$$
\par
We are in position to give 
the notion of weak solution.
\begin{defi}\label{defnomu}
{\rm 
Let $u_\circ\in L^2(\O)$. A field
$u:(0,\infty)\times \O\to \R^N$ 
 is said a weak solution of
system {\rm \eqref{PF}} if, for all $T\in (0,\infty)$,  $u\in W^p(\O_T)$,
$$\ba{ll}\dy \vs1 
\int_{0}^{T} [<u_\tau,\psi>_{V^p_\tau(\O)}+\left(|\nabla u|^{p-2}\nabla u,\nabla \psi\right)]dt =0\,,  \hskip0.1cm \ \forall \psi\in
X^p(\O_T)\,, \\\dy
\lim_{t\to 0^+}\|u(t)-u_\circ\|_2=0\,.
\ea
$$}
\end{defi}
In order to prove Theorem \ref{mmtp} we introduce the approximating systems
\be\label{PFl}\begin{array}{ll}\dy\vs1
u_t-\nabla\cdot\left(\left(\mu+|\nabla
u|^2\right)^\frac{p-2}{2}\nabla u\right)= 0\,,&\hskip-0.2cm\textrm{ in
}(0,\infty)\times\O,
\\\dy \vs1\hskip4cmu(t,x)=0\,,&\hskip-0.2cm\textrm{ on }
(0,\infty)\times\po,\;\\\dy\hskip4cm u(0,x)=u_\circ(x),&\hskip-0.2cm\mbox{ in
}\O,\end{array}\ee with $\mu\in (0,1)$, and

\be\label{PFepv}\begin{array}{ll}\vs1\dy v_t-\nu\Delta
v-\nabla\cdot\left(\left(\mu+|\nabla
v|^2\right)^\frac{p-2}{2}\nabla v\right)= 0\,,&\hskip-0.2cm\textrm{ in
}(0,\infty)\times\O,
\\\dy\vs1\hskip5cm v(t,x)=0\,,&\hskip-0.2cm \textrm{ on }
(0,\infty)\times\po,\;\\\dy\hskip5cm v(0,x)=u_\circ(x),&\hskip-0.2cm\mbox{ in
}\O\,.\end{array}\ee 
 \begin{defi}\label{wmu} {\rm Let $\mu\in (0,1)$. Let $u_{\circ}\in
L^2(\O)$. A field
$u:(0,\infty)\times \O\to \R^N$
 is said a weak solution of
system {\rm \eqref{PFl}} if, for all $T\in (0,\infty)$,  $u\in W^p(\O_T)$, $$\ba{ll}\dy \vs1 
\int_{0}^{T} [<u_\tau,\psi>_{V^p_\tau(\O)}+\left(a(\mu,u)\nabla u,\nabla\psi\right)]d\tau=0\,, \ \forall \psi\in
X^p(\O_T)\,,
\ea$$
$$\lim_{t\to
0^+}\|u(t)-u_\circ\|_2=0\,.$$}
\end{defi}
Set $q=q(t,x):=\max\{2,p(t,x)\}$. 
\begin{defi}\label{wmunu}
{\rm Let $\mu\in (0,1)$, $\nu>0$. Let $u_{\circ}\in L^2(\O)$. A field
$v:(0,\infty)\times \O\to \R^N$  is said a weak solution of
system {\rm \eqref{PFepv}} if,  for all $T\in (0,\infty)$, $v\in L^2(0,T; W_0^{1,2}(\O))\cap 
W^q(\O_T)$,
$$\ba{ll}\dy 
\int_{0}^{T}\hskip-0.2cm [<v_\tau,\psi>_{V^q_\tau(\O)}+\nu(\nabla v,\nabla\psi)+\left(a(\mu,v)\nabla v,\nabla\psi\right)]d\tau=0\,, \  
\ \forall \psi\!\in\!
X^q(\O_T)\,, 
\ea$$
 $$\lim_{t\to
0^+}\|v(t)-u_\circ\|_2=0\,.$$}
\end{defi}

System \eqref{PFl} is introduced in order to deal with a non-singular system. The introduction of system \eqref{PFepv} is connected with the idea of applying a duality technique to prove Theorem \ref{mmtp}. Let us explain in what sense. \par
In order to apply the duality, one aims at a reciprocity relation between the solution $u$ and the solution of a local adjoint problem. Both should be in the sets of solutions, that is the space $W^p(\O_T)$ and $W^{\ov p}(\O_t)$ respectively, $\ov p=\ov p(s,x):=p(t-s,x)$. The ``natural'' adjoint of system \eqref{PFl} should be deduced, on $(0,t)$, from the following
\be\label{AD01}\begin{array}{ll}\dy\vs1
  \psi_s- \nabla\cdot (B(u)(s,x)\nabla
  \psi)=
  0\,,&\hskip-0.2cm\textrm{ in }(0,t)\times\O,
  \end{array}\ee
where,  for all $t>0$, 
 we define for a.e. in $s\in (0,t)$, 
\be\label{beeta}
\ba{ll}\dy\vs1
B(u)(s,x)=(B(u)(s,x))_{i\alpha j\beta}\!:=
\delta_{ij}\,\delta_{\alpha\beta}(\mu+\!|\nabla u(t\!-\!s,x)|^2)^\frac
{p(t-s,x)-2}{2}\,,\ea
\ee 
with $\mu\in (0,1)$. This system is not suitable, as $\psi\not\in X^{\ov p}(\O_t)$. An approximation of the adjoint which ensures the membership to $X^{\ov p}(\O_t)$ is  
\be\label{AD001}\begin{array}{ll}\dy\vs1
  \theta_s- \nabla\cdot (B(u)(s,x)\nabla
  \theta)-\ve\nabla\cdot(|\nabla\theta|^{\ov p-2}\nabla\theta)=
  0\,,&\hskip-0.2cm\textrm{ in }(0,t)\times\O,
  \end{array}\ee
where $\ve>0$. In this case $\theta\in X^{\ov p}(\O_t)$ but unfortunately $\theta_s\not\in (X^{\ov p}(\O_t))'$. In order to overcome this {\it impasse}, we are led to modify the same system \eqref{PFl}, approximating, in turn, with system \eqref{PFepv}. Now, the solution of \eqref{PFepv} belongs to $W^q(\O_T)$, with $q=q(t,x):=\max\{2,p(t,x)\}$. Then, the same arguments that led us to system \eqref{AD001} now lead us to the following one as approximation of the adjoint of \eqref{PFepv}
$$\begin{array}{ll}\dy\vs1
  \vp_s-\nu\Delta\vp- \nabla\cdot (B(v)(s,x)\nabla
  \vp)-\ve\nabla\cdot(|\nabla\vp|^{\ov q-2}\nabla\vp)=
  0\,, \textrm{ in }(0,t)\times\O, \ea$$
 where  $ \ov q=\ov q(s,x):=q(t-s,x) $. The solution of such system are showed to 
  belong to the right space $W^{\ov q}(\O_t)$, so that the following reciprocity relation can be obtained
\be\label{a9}\ba{ll}\dy\vs1 (v(t),\vp_\circ)=(v(0),\vp^{\ve}(t))\!-\ve
\int_0^t(|\nabla \hat \vp^{\ve}(\tau)|^{q(\tau,x)-2}\nabla \hat
\vp^{\ve}(\tau), \nabla v)d\tau. \hfill\dy
\ea\ee
This  relation enables to get global $L^{r_0}$ estimates for solutions of system \eqref{PFepv}. Since this system depends on two parameters, $\mu>0$ and $\nu>0$, and in the limit as $\mu\to 0$ and $\nu\to 0$ leads to system \eqref{PF}, the last step is to perform such limit procedures. \vskip0.2cm 
The proof of Theorem \ref{mmtpx} and its corollary is based, as in \cite{HV, Ver}, on the study of the following differential inequality
\be\label{per1}\frac1{r_0} \frac{d}{d\tau}\|u(\tau)\|_{r_0}^{r_0}+\int_\O |\nabla u(\tau,x)|^{p(\tau,x)}|u(\tau,x)|^{r_0-2}dx\leq 0, \ \mbox{a.e. in } (0,T)\,,\ee
 together with an iterative argument. In order to get \eqref{per1} firstly we approximate the data, and hence the solution, with a sequence of initial data in $C_0^\infty(\O)$. 
Since $u_\circ \in C_0^\infty(\O)$ implies
$u_\circ \in L^\infty(\O)$, it follows from Theorem \ref{mmtp}
that $u(t)\in L^\infty(\O)$, for any $t\geq 0$, and $u \in
C([0,T];L^{\bar r}(\Omega))$ for any $\bar r <\infty$. Then one readily
recognizes that $|u|^{r_0-2} u\in X^p(\O_T)$, for any $r_0\in
[2,\infty)$, so that $u|u|^{r_0-2}$ can be used as test function in the weak formulation of
\eqref{PF}, and one readily gets \eqref{per1}. If $u_\circ$ is just in $L^{r_0}(\O)$, denoting by $\{u_{\circ}^n\}$ a sequence of functions in $C_0^\infty(\O)$ strongly converging to $u_\circ$ in $L^{r_0}(\O)$ and such that $\|u_\circ^n\|_{r_0}\leq c\|u_\circ\|_{r_0}$, and denoting by $\{u^n\}$ the corresponding sequence of weak solutions to problem \eqref{PF}, for the sequence $\{u^n\}$ one finds  
  inequality \eqref{mmtp2} with $u^n$ in place of $u$.     
        Then, as the estimate is uniform in $n\in \N$, we can take a subsequence converging to a limit function $u(t)$, that satisfies the same estimate. Finally, that the limit $u$ is solution of \eqref{PF} corresponding to the initial data $u_\circ\in L^{r_0}(\O)$ is standard. \par 
 We like to point out that a different approach could be the study of suitable integral inequalities, as made for instance in \cite{MP1,MP2} in order to prove the regularizing effect in the scalar case with a constant exponent $p$. However, due to the vectorial character of our problem and to the variable exponent $p(t,x)$, an extension of the quoted technique to our setting seems to be not straightforward.

\end{document}